\documentclass[reqno,11pt]{amsart} 
\usepackage{amsmath}
\usepackage{amssymb}
\usepackage{amsthm}

\newcommand\NoBlackBoxes{\global\overfullrule0pt}
\NoBlackBoxes
\theoremstyle{plain} 

\def\4{\kern1pt}

\def\6{\vphantom0}

\def\8{\kern-10pt}
\def\7#1{_{(#1)}}

\begin{document}

\def\ffrac#1#2{\raise.5pt\hbox{\small$\4\displaystyle\frac{\,#1\,}{\,#2\,}\4$}}
\def\ovln#1{\,{\overline{\!#1}}}
\def\ve{\varepsilon}
\def\kar{\beta_r}

\title{NON-UNIFORM BOUNDS IN THE POISSON  APPROXIMATION \\ 
WITH APPLICATIONS TO INFORMATIONAL DISTANCES. I \\
}

\author{S. G. Bobkov$^{1}$}
\thanks{1) 
School of Mathematics, University of Minnesota, USA;
research was partially supported by SFB 1283, Humboldt Foundation, 
and NSF grant}
\address
{Sergey G. Bobkov \newline
School of Mathematics, University of Minnesota  \newline 
127 Vincent Hall, 206 Church St. S.E., Minneapolis, MN 55455 USA}
\smallskip
\email{bobkov@math.umn.edu}

\author{G. P. Chistyakov$^{2}$}
\thanks{2) 
Faculty of Mathematics, University of Bielefeld, Germany;
research was partially supported by SFB 1283}
\address
{Gennadiy P. Chistyakov\newline
Fakult\"at f\"ur Mathematik, Universit\"at Bielefeld\newline
Postfach 100131, 33501 Bielefeld, Germany}
\smallskip
\email{goetze@mathematik.uni-bielefeld.de}

\author{F. G\"otze$^{2}$}
\address
{Friedrich G\"otze\newline
Fakult\"at f\"ur Mathematik, Universit\"at Bielefeld\newline
Postfach 100131, 33501 Bielefeld, Germany}
\email{goetze@mathematik.uni-bielefeld.de}


\begin{abstract}
We explore asymptotically optimal bounds for deviations of Bernoulli
convolutions from the Poisson limit in terms of the Shannon relative entropy 
and the Pearson $\chi^2$-distance. The results are based 
on proper non-uniform estimates for densities. This part deals
with the so-called non-degenerate case.
\end{abstract}

\maketitle
\markboth{S. G. Bobkov, G. P. Chistyakov and F. G\"otze}{
Relative entropy and $\chi^2$ divergence from the Poisson law}

\def\theequation{\thesection.\arabic{equation}}
\def\E{{\mathbb E}}
\def\R{{\mathbb R}}
\def\C{{\mathbb C}}
\def\P{{\mathbb P}}
\def\Z{{\mathbb Z}}
\def\H{{\rm H}}
\def\Im{{\rm Im}}
\def\Tr{{\rm Tr}}

\def\k{{\kappa}}
\def\M{{\cal M}}
\def\Var{{\rm Var}}
\def\Ent{{\rm Ent}}
\def\O{{\rm Osc}_\mu}

\def\ep{\varepsilon}
\def\phi{\varphi}
\def\F{{\cal F}}
\def\L{{\cal L}}

\def\be{\begin{equation}}
\def\en{\end{equation}}
\def\bee{\begin{eqnarray*}}
\def\ene{\end{eqnarray*}}

\section{{\bf Introduction}}
\setcounter{equation}{0}

\vskip2mm
\noindent
Let $X_1,\dots,X_n$ be independent Bernoulli random variables taking 
the two values, $1$ (interpreted as a success) and $0$ (as a failure) 
with respective probabilities $p_j$ and $q_j = 1-p_j$. The total number 
of successes $W = X_1 + \dots + X_n$ takes values $k = 0,1,\dots,n$ 
with probabilities
\be
\P\{W = k\} = 
\sum p_1^{\ep_1} q_1^{1-\ep_1} \dots p_n^{\ep_n} q_n^{1-\ep_n},
\en
where the summation runs over all 0-1 sequences $\ep_1,\dots,\ep_n$ 
such that $\ep_1 + \dots + \ep_n = k$. Although this expression is difficult 
to determine in case of  arbitrary $p_j$ and large $n$, it can be well 
approximated by the Poisson probabilities under quite general assumptions. 
Putting 
$$
\lambda = p_1 + \dots + p_n,
$$
let $Z$ be a Poisson random variable with parameter $\lambda>0$
(for short, $Z \sim P_\lambda$), i.e.,
$$
v_k = \P\{Z = k\} = \frac{\lambda^k}{k!}\,e^{-\lambda}, \qquad k = 0,1,\dots
$$
It is well-known for a long time that, if $\max_{j \leq n} p_j$ 
is small, the distribution $P_\lambda$ approximates the distribution 
$P_W$ of $W$, which may be quantified by means of the total variation 
distance
\bee
d(W,Z) 
 & = & 
\|P_W - P_\lambda\|_{\rm TV} \\
 & = & 
2\,\sup_{A \subset \Z}\, |\,\P\{W \in A\} - \P\{Z \in A\}|
 \, = \,
\sum_{k=0}^\infty\, |w_k - v_k|,
\ene
where $w_k = \P\{W = k\}$. In particular, 
based on Stein-Chen's method, there is the following two-sided 
bound due to Barbour and Hall involving the functional
$$
\lambda_2 = p_1^2 + \dots + p_n^2.
$$

\vskip5mm
{\bf Theorem 1.1} \cite{B-H}. {\sl One has
\be
\frac{1}{32} \min(1,1/\lambda)\,\lambda_2 \, \leq \,
\frac{1}{2}\, d(W,Z) \, \leq \,
\frac{1-e^{-\lambda}}{\lambda}\,\lambda_2.
\en
}

\vskip2mm
Here, the parameter $\lambda_2$, or more precisely -- the ratio 
$\lambda_2/\lambda$ (for $\lambda$ bounded away from zero), plays 
a similar role as the Lyapunov ratio $L_3$ in the central limit theorem.

In the i.i.d. case with $p_j = \lambda/n$ and  fixed $\lambda>0$, 
both sides of (1.2) are of the same order $1/n$. In the case 
$\lambda \leq 1$, 
the upper bound in (1.2) is sharp also in the sense that the second inequality 
becomes an equality for $p_1 = \lambda$, $p_j = 0$ ($2 \leq j \leq n$).

Theorem 1.1 refined many previous results in this direction, 
starting from bounds for the i.i.d. case by Prokhorov \cite{P} and bounds 
for the general case by Le Cam \cite{LC1}. In particular, Le Cam obtained 
the upper bound
\be
d(W,Z) \leq 2\lambda_2.
\en
For large $\lambda$ Kerstan \cite{K} and respectively Chen \cite{Ch} improved 
these bounds to
$$
d(W,Z) \, \leq \, \frac{2.1}{\lambda}\,\lambda_2
\quad {\rm if} \ \ \max_{j \leq n} p_j \leq \frac{1}{4}, \qquad
{\rm respectively} \ \ 
d(W,Z) \, \leq \, \frac{10}{\lambda}\,\lambda_2.
$$
See also \cite{H-LC}, \cite{Ve}, \cite{Se}, \cite{R}, \cite{Ro}, \cite{B-H-J} 
and the references therein. A certain refinement of the lower bound in (1.2) 
was obtained in Sason \cite{S1}.

While (1.2) provides a sharp estimate for the  total variation distance, 
one may wonder whether or not similar approximation bounds hold for
the stronger informational distances. As a first interesting example, one
may consider the relative entropy
$$
D(W||Z) = \sum_{k=0}^\infty w_k \log \frac{w_k}{v_k},
$$
often called the Kullback-Leibler distance, or an informational divergence 
of $P_W$ from $P_\lambda$.  It dominates the total variation distance
in view of the Pinsker inequality
$
D(W||Z) \geq \frac{1}{2}\,d(W,Z)^2.
$
In this context, lower and upper bounds for the relative entropy 
were studied by Harremo\"es \cite{H2001}, \cite{H2003}, and
Harremo\"es and Ruzankin \cite{H-R}. In particular, 
in the i.i.d. case $p_j = p$, it was shown in \cite{H-R} that
\bee
\frac{-\log(1-p) - p}{2} - \frac{14 p^2}{n\,(1-p)^3} 
 & \leq & 
D(W||Z) \\
 & \leq & 
\frac{-\log(1-p) - p}{2} - \frac{(1+p)\, p^2}{4n\,(1-p)^3}.
\ene
If $p = \lambda/n$ with a fixed (or just bounded) value of $\lambda$, 
these estimates provide the rate of Poisson approximation
\be
D(W||Z) \, = \, \frac{\lambda^2}{4n^2} + O(1/n^3) \quad {\rm as} \ \ 
n \rightarrow \infty.
\en

The general non-i.i.d. scenario (with not necessarily equal probabilities $p_j$) 
has been partially studied as well. A simple upper estimate
$D(W||Z) \leq \lambda_2$, analogous to Le Cam's bound (1.3),
may be found in \cite{H2001}, cf. also Johnson \cite{J}. It is however
not so sharp as (1.4). A tighter upper bound
\be
D(W||Z) \, \leq \, \frac{1}{\lambda}\,\sum_{j=1}^n \frac{p_j^3}{1 - p_j}
\en
was later derived by Kontoyiannis, Harremo\"es and Johnson \cite{K-H-J}.
If $p_j = \lambda/n$ with $\lambda \leq n/2$, it yields
$D(W||Z) \leq 2\lambda^2/n^2$ reflecting a correct decay with respect to $n$ 
up to a constant, according to (1.4). Nevertheless, in the general case, 
Pinsker's inequality and the bounds (1.2)-(1.3) suggest that a further 
sharpening such as 
\be
D(W||Z) \leq A_\lambda \lambda_2^2
\en
might be possible by involving $\lambda_2$ rather than the functional 
$\lambda_3 = p_1^3 + \dots + p_n^3$. To compare the two quantities, note that, 
by Cauchy's inequality, $\lambda_2^2 \leq \lambda \lambda_3$. Hence, the
inequality (1.6) would be sharper compared  to (1.5) modulo 
a $\lambda$-dependent factor. An upper bound such as (1.6) may also be 
inspired by the lower bound
\be
D(W||Z) \geq \frac{1}{4}\, \Big(\frac{\lambda_2}{\lambda}\Big)^2
\en
recently derived by Harremo\"es, Johnson and Kontoyiannis \cite{H-J-K}.
It is consistent with (1.4) and also shows that the constant $1/4$
is best possible.

As it turns out, (1.6) does hold in the so-called non-degenerate situation,
and in essence, (1.7) may be reversed
(we say that the range of $(\lambda,\lambda_2)$ is non-degenerate, if 
$\lambda_2 \leq \kappa \lambda$ with $\kappa \in (0,1)$,
or if $\lambda \leq \lambda_0$, and implicitly mean that the resulting
inequalities may contain $\kappa$ or $\lambda_0$ as fixed parameters).
Moreover, one can further sharpen (1.6) by replacing the relative 
entropy with the Pearson $\chi^2$-distance, as well as with other 
R\'enyi/Tsallis distances. To avoid technical complications,
let us restrict ourselves to the $\chi^2$-divergence which is given by
$$
\chi^2(W,Z) = \sum_{k=0}^\infty 
\frac{(w_k - v_k)^2}{v_k}.
$$
It is a divergence type quantity which dominates the relative entropy:
$\chi^2(W,Z) \geq D(W||Z)$.
For a general theory of informational distances, we refer interested readers
to the recent review by van Erven and Harremo\"es \cite{E-H}; an additional 
material may be found in the books \cite{LC2}, \cite{L-V}, \cite{Va}, \cite{J}.
Here, we reverse the bound (1.7) and prove:

\vskip5mm
{\bf Theorem 1.2.} {\sl If $\lambda_2 \leq \lambda/2$, then
with some absolute constant $c$ we have
\be
D(W||Z) \, \leq \, \chi^2(W,Z) \, \leq \, c\,
\Big(\frac{\lambda_2}{\lambda}\Big)^2.
\en 
}

\vskip2mm
The condition $\lambda_2 \leq \lambda/2$ is readily fulfilled as long as 
all $p_j \leq 1/2$ (note that, if $\lambda \leq 1/2$, then necessarily 
$p_j \leq 1/2$ and then $\lambda_2 \leq \lambda/2$). Similar bounds as in (1.8) remain to hold under 
the weaker assumption $\lambda_2 \leq \kappa \lambda$ with a constant 
$c = c_\kappa$ depending on $\kappa \in (0,1)$, cf. Proposition 6.2 below. 
This assumption may actually be replaced with the requirement that $\lambda$ 
is bounded. More precisely, in the second part of the paper it will be 
shown that without any restriction, up to some universal factors, we have
\bee
D(W||Z) \sim
\Big(\frac{\lambda_2}{\lambda}\Big)^2\,(1 + \log F), \qquad
\chi^2(W,Z) \sim
\Big(\frac{\lambda_2}{\lambda}\Big)^2\,\sqrt{F},
\ene
where
$$
F = \frac{\max(1,\lambda)}{\max(1,\lambda-\lambda_2)}.
$$
This shows that in general the bound (1.7)  cannot be reversed.

For the study of the asymptotic behavior of $D$ and $\chi^2$ in terms of 
$\lambda$ and $\lambda_2$, we derive new bounds for the difference
between densities of $W$ and $Z$, that is, for
$$
\Delta_k = w_k - v_k = \P\{W=k\} - \P\{Z=k\}.
$$
To this aim, one has to consider different zones of $\lambda$'s, 
distinguishing between ``small'' and ``large'' values. The case 
$\lambda \leq \frac{1}{2}$ can be handled directly leading to 
the non-uniform density bound
$$
|\Delta_k| \, \leq \, 2 \lambda_2\,\P\{k-2 \leq Z \leq k\}.
$$
It easily yields sharp upper bounds for all above distances 
as in Theorems 1.1-1.2 in the case of small $\lambda$, at least up to
numerical factors (cf. Proposition 3.3 and 3.4).
To treat larger values of $\lambda$, a more sophisticated analysis in 
the complex plane is involved -- using the closeness of the generating 
functions associated with the sequences $w_k$ and $v_k$. 
In particular, the following statement may be of independent interest.

\vskip5mm
{\bf Theorem 1.3.} {\sl For all integer $k \geq 0$, 
\be
|\Delta_0|\leq 3\lambda_2\,e^{-\lambda},
\quad |\Delta_k| \, \leq \, 3 \lambda_2\,\,\, (k\geq 1).
\en
Moreover, putting 
$\rho = (\lambda-\lambda_2)\,\min\{\frac{k}{\lambda},\frac{\lambda}{k}\}$,
$k = 1,2,\dots$, we have
\begin{eqnarray}
|\Delta_k|
 & \le & 
7 \sqrt{k}\,\Big(\frac{k-\lambda}{\lambda}\Big)^2\,
\lambda_2\min\big\{1,\rho^{-1/2}\big\}\ \P\{Z=k\}\nonumber \\
 & & + \,
21\,k^{3/2}\,
\frac{\lambda_2}{\lambda}\min\big\{1,\rho^{-3/2}\big\}\ \P\{Z=k\}.
\end{eqnarray}
}

\vskip2mm
Let us clarify the meaning of the last bound, assuming that 
$\lambda_2 \leq \kappa \lambda$ with some constant $\kappa \in (0,1)$. 
If $k \leq 2\lambda$ and $\lambda \geq 1/2$, then with some 
$c = c_\kappa > 0$, it gives
$$
|\Delta_k| \, \leq \,
c\, \Big(\frac{(k-\lambda)^2}{\lambda} + 1\Big)\, 
\frac{\lambda_2}{\lambda}\, \P\{Z=k\},
$$
while for $k \geq \lambda \geq 1/2$, we also have
$$
|\Delta_k| \, \leq \,
c\, \Big(\frac{k}{\lambda}\Big)^3\,  \lambda_2 \, \P\{Z=k\}.
$$
Since $|k - \lambda|$ is of order at most $\sqrt{\lambda}$ on a sufficiently 
large part of $\Z$ measured by $P_\lambda$, these non-uniform bounds 
explain the possibility of upper bounds in Theorem 1.2.

The paper is organized as follows. First we describe several general
bounds on the probability function of the Poisson law (Section 2).
In Sections 3, we consider the deviations $\Delta_k$ and prove Theorem 1.2 
in case $\lambda \leq 1/2$. Sections 4-5 are devoted to non-uniform bounds 
and the proof of Theorem 1.3, which is used to complete the proof of 
Theorem 1.2 for $\lambda \geq 1/2$. Uniform bounds for large $\lambda$ 
are discussed in Section 7. There we shall demonstrate that in a typical 
situation, when the ratio $\lambda_2/\lambda$ is small, the Poisson 
approximation considerably improves the rate of normal approximation 
described by the Berry-Esseen bound in the central limit theorem.

\vskip7mm
\section{{\bf Gaussian Type Bounds on Poisson Probabilities}}
\setcounter{equation}{0}

\vskip2mm
\noindent
When bounding the Poisson probabilities
$$
v_k = f(k) = \P\{Z = k\} = \frac{\lambda^k}{k!}\,e^{-\lambda}, \qquad 
k = 0,1,\dots,
$$
with a fixed parameter $\lambda>0$, it is convenient to use the well-known 
Stirling-type two-sided bound:
\be
\sqrt{2\pi}\,k^{k + \frac{1}{2}}\,e^{-k} \leq k! \leq \, 
e\,k^{k + \frac{1}{2}}\,e^{-k} \qquad (k \geq 1).
\en
In particular, it implies the following Gaussian type estimates.

\vskip5mm
{\bf Lemma 2.1.} {\sl For all $k \geq 1$,
\be
f(k) \, \leq \, \frac{1}{\sqrt{2\pi k}}.
\en
Moreover, if $1 \leq k \leq 2\lambda$, then
\be
\frac{1}{e\sqrt{k}}\,e^{-\frac{(k-\lambda)^2}{\lambda}} \, \leq \, 
f(k) \, \leq \, \frac{1}{\sqrt{2\pi k}}\,e^{-\frac{(k-\lambda)^2}{3\lambda}}.
\en
Here, the lower bound may be improved in the region $k \geq \lambda$ as
\be
f(k) \, \geq \, \frac{1}{e\sqrt{k}}\,e^{-\frac{(k-\lambda)^2}{2\lambda}}.
\en
}

{\bf Proof.} Applying the lower estimate in (2.1), we get
\begin{eqnarray}
f(k) 
 & \le &
\frac{1}{\sqrt{2\pi k}}\, e^{k - \lambda}\, \Big(\frac{\lambda}{k}\Big)^k \\
 & = &
\frac{1}{\sqrt{2\pi k}}\ e^{k -\lambda} \ 
e^{-k\log(1+\frac{k-\lambda}{\lambda})}  \, = \,
\frac{1}{\sqrt{2\pi k}}\ e^{\lambda h(\theta)}, \qquad 
\theta = \frac{k-\lambda}{\lambda} \nonumber,
\end{eqnarray}
where
$$
h(\theta) = \theta - (1 + \theta)\log(1 + \theta).
$$
The function $h(\theta)$ is concave on the half-axis $\theta \geq -1$, with
$h(0) = h'(0) = 0$. Hence, $h(\theta) \leq 0$
for all $\theta$, thus proving the first assertion (2.2).

Assuming that $1 \leq k \leq 2\lambda$ (with $\lambda \geq \frac{1}{2}$), 
we necessarily have $|\theta| \leq 1$. In this interval, consider 
the function $T_c(\theta) = h(\theta) + c\theta^2$ with parameter $\frac{1}{4}<c\le 1$.
The second derivative 
$$
T_c''(\theta) = - \frac{1}{1 + \theta} + 2c \qquad (-1 < \theta \leq 1)
$$
is vanishing at the point $\theta_0=\frac{1}{2c}-1$, while $T_c''(-1) = -\infty$. This means that $T_c$
is concave on $[-1,\theta_0]$ and convex on $[\theta_0,1]$.
Since also $T_c(0) = T_c'(0) = 0$, we have $T_c(\theta)\le 0$ for all $\theta\in[-1,1]$, if and only if this inequality is fulfilled at $\theta=1$.
But $T_c(1)=1-2\log 2+c$, so the optimal value is
$c=2\log 2-1=0.387...>\frac{1}{3}$. Hence,
$h(\theta) \leq -\frac{1}{3}\,\theta^2$, and we arrive at the upper bound in (2.3).

Similarly, applying the upper estimate in (2.1), we get
$$
f(k) \, \geq \,
\frac{1}{e\sqrt{k}}\ e^{k - \lambda}\, \Big(\frac{\lambda}{k}\Big)^k 
 \, = \, \frac{1}{e\sqrt{k}}\ e^{\lambda h(\theta)}, \qquad 
\theta = \frac{k-\lambda}{\lambda}.
$$
Choosing $c=1$, consider the function
$T(\theta) = h(\theta) + \theta^2$ in the interval $|\theta| \leq 1$. 
Since $T''(-\frac{1}{2}) = 0$, it is concave on $[-1,-\frac{1}{2}]$ 
and is convex on $[-\frac{1}{2},1]$. Since $T(0) = T'(0) = 0$ and 
$T(-1) = 0$, this means that $\theta = 0$ is the point of minimum of $T$. 
Therefore, $T(\theta) \geq 0$, that is, 
$h(\theta) \geq -\theta^2$ for all $\theta \in [-1,1]$, giving
the lower bound in (2.3).

Finally, to get the refinement (2.4) in the region $k \geq \lambda$, 
consider the function
$T(\theta) = h(\theta) + \frac{1}{2}\,\theta^2$ for $\theta \geq 0$. 
Since $T(0) = 0$ and $T'(\theta) = \theta - \log(1 + \theta) \geq 0$,
this function is increasing. Therefore, $T(\theta) \geq 0$, that is, 
$h(\theta) \geq -\frac{1}{2}\,\theta^2$ for all $\theta \geq 0$.
\qed

\vskip7mm
\section{{\bf \large Elementary Upper Bounds}}
\setcounter{equation}{0}

\vskip2mm
\noindent
We keep the same notations as before; in particular,
$$
\P\{Z = k\} = \frac{\lambda^k}{k!}\,e^{-\lambda}, \qquad k = 0,1,\dots,
$$
while
$$
\P\{W = k\} \, = \, \sum p_1^{\ep_1}\, 
(1-p_1)^{1-\ep_1} \dots p_n^{\ep_n}\, (1-p_n)^{1-\ep_n}
$$
with summation over all 0-1 sequences $\ep = (\ep_1,\dots,\ep_n)$
such that $\ep_1 + \dots + \ep_n = k$. Clearly, $\P\{W = k\} = 0$
for $k > n$. To eliminate this condition,
one may always assume that $n$ is arbitrary, by extending the sequence
$(X_1,\dots,X_n)$ to $(X_1,\dots,X_k)$ in case $n < k$ with
$p_{n+1} = \dots = p_k = 0$. This does not change the value of $W$.

First, let us consider the probability that $W$ equals 
$k=0$.

\vskip5mm
{\bf Lemma 3.1.} {\sl If $\max_j p_j \leq \frac{1}{2}$, then
\bee
0 \, \leq \, \P\{Z=0\} - \P\{W=0\} 
 & \leq &
\lambda_2\,e^{-\lambda}. 
\ene
}

{\bf Proof.} Expanding the function $p \rightarrow -\log(1-p)$ near zero 
according to the Taylor formula as in the previous section, write
\be
\P\{W=0\} =
\prod_{j=1}^n (1-p_j) = e^{-\lambda - S}, \qquad
S = \sum_{s=2}^\infty \ \frac{1}{s}\,\lambda_s,
\quad \lambda_s=p_1^s+\dots+p_n^s.
\en
Using 
$\lambda_s \leq (\max_j p_j)^{s-2}\,\lambda_2 \leq 2^{-(s-2)}\,\lambda_2$ 
for $s \geq 2$, we have
\be
S \, \leq \, \lambda_2 \sum_{s=2}^\infty \frac{2^{-(s-2)}}{s} \, = \,
(4\log 2 - 2)\,\lambda_2 \, \leq \,\lambda_2.
\en
Hence
$$
\P\{Z=0\} - \P\{W=0\} = e^{-\lambda}\,(1 - e^{-S}) \leq e^{-\lambda}\,S.
$$

\qed

\vskip2mm
Note that the condition of Lemma 3.1 is fulfilled automatically, if
$\lambda \leq 1/2$. In that case, the upper bounds of the lemma may easily
be reversed up to numerical factors, for example, in the form
\bee
\P\{Z=0\} - \P\{W=0\} 
 & \geq &
0.47\,\lambda_2\,e^{-\lambda}, \\
\P\{W=1\} - \P\{Z=1\} 
 & \geq &
0.42\,\lambda_2\,e^{-\lambda}.
\ene
Moreover, if $\lambda \leq 1/8$, then also
$$
\P\{Z=2\} - \P\{W=2\} \, \geq \,
\frac{17}{49}\,\lambda_2\,e^{-\lambda}.
$$
Here, the value $k=2$ turns out to be most essential for obtaining
lower bounds, since it immediately yields $d(W,Z) \geq c\lambda_2$ and
$D(W||Z) \geq c\,(\frac{\lambda_2}{\lambda})^2$
with some absolute constant $c>0$.

Returning to upper bounds, recall the notation
$\Delta_k=\P\{W=k\}-\P\{Z=k\}$.
In order to involve the values $k \geq 1$, we need the following:

\vskip5mm
{\bf Lemma 3.2.} {\sl If $\max_j p_j \leq 1/2$, then 
\be
|\Delta_1|\le\lambda_2(\lambda+e-1)e^{-\lambda}.
\en
Moreover, for any $k\ge 2$,
$$
\big|\P\{W=k\} - \P\{Z=k\}\big| \, \leq \, \lambda_2\,
\Big(\frac{\lambda^k}{k!} + 
\frac{e^\lambda - 1}{\lambda}\,\frac{\lambda^{k-1}}{(k-1)!} + 
\frac{\lambda^{k-2}}{(k-2)!}\Big)\,e^{-\lambda}.
$$
}

\vskip2mm
{\bf Proof.} 
Denote by $I$ the collection of all tuples $\varepsilon=(\varepsilon_1,\dots,\varepsilon_n)$
with integer components $\varepsilon_i\ge 0$ such $\varepsilon_1+\dots+\varepsilon_n=k$, and let $J=\{\varepsilon\in I:\max_i\varepsilon_i\le 1\}$.
Representing the Poisson random variable $Z \sim P_\lambda$ as
$Z = Z_1 + \dots + Z_n$ with independent summands
$Z_j \sim P_{p_j}$, we have that, for any $k = 0,1,\dots$, 
$$
\P\{Z=k\} = e^{-\lambda} \sum_{\ep\in I}
\frac{p_1^{\ep_1} \dots p_n^{\ep_n}}{\ep_1!\dots\ep_n!}.
$$
Hence, we may start with the formula
\be
\P\{Z=k\} - \P\{W=k\} \, = \, e^{-\lambda}
\sum_{\ep\in I}
\frac{1}{\ep_1!\dots\ep_n!}\ U_\ep \ \ 
- \sum_{\ep \in J} U_\ep V_\ep,
\en
where
$$
U_\ep =  p_1^{\ep_1}\dots p_n^{\ep_n}, \qquad
V_\ep = (1-p_1)^{1-\ep_1} \dots (1-p_n)^{1-\ep_n}.
$$

For a 0-1 sequence $\ep = (\ep_1,\dots,\ep_n)\in J$, put
$$
L_\ep = \ep_1 p_1 + \dots + \ep_n p_n.
$$
By the Taylor formula once more, 
$$
V_\ep^{-1} \, = \, e^{S_\ep}, \qquad
S_\ep = \sum_{s=1}^\infty \ \frac{1}{s}\, \sum_{j=1}^n (1-\ep_j)\, p_j^s.
$$
Similarly to (3.1)-(3.2), we have
$$
S_\ep \, = \, \lambda - L_\ep + 
\sum_{s=2}^\infty \ \frac{1}{s}\, \sum_{j=1}^n\, (1-\ep_j)\, p_j^s \, = \, 
\lambda - L_\ep + \theta \lambda_2, \qquad 0 \leq \theta \le 1.
$$
Therefore,
$$
e^\lambda\, V_\ep = e^{L_\ep - \theta \lambda_2}
\geq 1 + (L_\ep - \theta \lambda_2) \geq 1 + L_\ep - \lambda_2.
$$
Moreover, since $L_\ep \leq \min{(\lambda,k)}$, we have
$\frac{e^{L_\ep} - 1}{L_\ep} \leq \frac{e^{\min(\lambda,k)} - 1}{\min(\lambda,k)} \equiv c_{k,\lambda}$, 
which in turn implies
$e^\lambda\, V_\ep \leq e^{L_\ep} \leq 1 + c_{k,\lambda}\, L_\ep$. 
The two bounds give
$
L_\ep - \lambda_2 \leq e^\lambda\, V_\ep - 1 \leq c_{k,\lambda}\, L_\ep,
$
so that
$$
\big|U_\ep - e^\lambda\,U_\ep V_\ep\big| \, \leq \, \lambda_2\,U_\ep +
c_{k,\lambda}\, U_\ep L_\ep.
$$

Next, applying the multinomial formula, we have
$$
\sum_{\ep\in J}
U_\ep \ \leq 
\sum_{\ep\in I}
\frac{p_1^{\ep_1}\dots p_n^{\ep_n}}{\ep_1!\dots\ep_n!} \, = \,
\frac{\lambda^k}{k!}
$$
and
\bee
\sum_{\ep\in J} U_\ep\,L_\ep
 & = &
\sum_{i=1}^n \ \sum_{\ep\in J}
\ep_i \ p_1^{\ep_1}\dots p_{i-1}^{\ep_{i-1}}\, p_i^{\ep_i + 1}\,
p_{i+1}^{\ep_{i+1}} \dots p_n^{\ep_n} \\
 & = & 
\sum_{i=1}^n p_i^2 
\sum_{\ep\in J, \, \ep_i=1}
p_1^{\ep_1}\dots p_{i-1}^{\ep_{i-1}}\,
p_{i+1}^{\ep_{i+1}} \dots p_n^{\ep_n} \\
 & \leq &
\sum_{i=1}^n p_i^2 \, 
\frac{1}{(k-1)!}\,(\lambda - p_i)^{k-1} \ \leq \ 
\lambda_2\, \frac{\lambda^{k-1}}{(k-1)!}. 
\ene
Thus,
\be
\sum_{\ep\in J} 
|U_\ep - e^{\lambda} \,U_\ep V_\ep|
 \, \leq \,
\lambda_2\,\Big(\frac{\lambda^k}{k!} + c_{k,\lambda}\, 
\frac{\lambda^{k-1}}{(k-1)!}\Big).
\en

The remaining terms participating in $\P(Z=k)$
correspond to the tuples $\ep\in I$ with $\max_i\ep_i\ge 2$, which is only possible for $k\ge 2$. In that case, restricting for definiteness to the constraint
$\ep_n\ge 2$, we have
\bee
\sum_{\ep\in I,\,\ep_n\ge 2}
\frac{p_1^{\ep_1}\dots p_n^{\ep_n}}{\ep_1!\dots\ep_n!}
 & = &
\sum_{m=2}^k \frac{p_n^m}{m!} 
\sum_{\ep_1 + \dots + \ep_{n-1} = k-m}
\frac{p_1^{\ep_1}\dots p_{n-1}^{\ep_{n-1}}}{\ep_1!\dots\ep_{n-1}!} \\
 & = &
\sum_{m=2}^k \frac{p_n^m}{m!} \ \frac{(\lambda - p_n)^{k-m}}{(k-m)!} \\ 
 & \le &
p_n^2 \sum_{m=2}^k \frac{p_n^{m-2}}{(m-2)!} \ \frac{(\lambda - p_n)^{k-m}}{(k-m)!}
 \ = \
p_n^2\,\frac{\lambda^{k-2}}{(k-2)!}.
\ene
Similarly, for any $i = 1,\dots,n$,
$$
\sum_{\ep\in I,\,\ep_i\ge 2}
\frac{p_1^{\ep_1}\dots p_n^{\ep_n}}{\ep_1!\dots\ep_n!} \, \leq \,
p_i^2\,\frac{\lambda^{k-2}}{(k-2)!},
$$
and summing over $i \leq n$, we then get
$$
\sum_{\ep\in I\,\max \ep_j\ge 2}
\frac{p_1^{\ep_1}\dots p_n^{\ep_n}}{\ep_1!\dots\ep_n!}
 \, \leq \,
\lambda_2\,\frac{\lambda^{k-2}}{(k-2)!}.
$$

\vskip5mm

It remains to combine this bound with the bound (3.5)
and apply both in (3.4). Then we finally obtain that
\be
|\Delta_k|\le \lambda_2\Big(\frac{\lambda^k}{k!}+
c_{k,\lambda}\frac{\lambda^{k-1}}{(k-1)!}+I_{\{k\ge 2\}}\frac{\lambda^{k-2}}{(k-2)!}\Big)\,e^{-\lambda}.
\en
If $k=1$, then $c_{1,\lambda}\le e-1$, and we arrive at the first inequality in (3.3). In the case $k\ge 2$, one may use $c_{k,\lambda}\le\frac{e^{\lambda}-1}{\lambda}$, and then we arrive at the second inequality of the lemma.
\qed

Note that when $\lambda\le\frac{1}{2}$, we also have $c_{1,\lambda}\le 2(\sqrt e-1)$, and then (3.3)
may be replaced with a slightly better bound
\be
|\Delta_1|\le\,\lambda_2(\lambda+2(\sqrt e-1))e^{-\lambda}.
\en

Combining Lemmas 3.1--3.2 (cf. (3.6)), we thus obtain the following non-uniform bound on the deviations of $\Delta_k$.

\vskip5mm
{\bf Proposition 3.3.} {\sl If $\max_j p_j\leq 1/2$, then, for all $k\geq 0$,
$$
|\Delta_k|\leq \frac{e^{\lambda}-1}{\lambda}\,\lambda_2\,\P\{k-2\le Z\le k\}.
$$
}

The estimates obtained so far are sufficient to establish Theorem 1.2 in the
case $\lambda \leq 1/2$. In fact, one may weaken the latter condition
to $\max_j p_j \leq 1/2$, as shown in the next statement.
To compare the lower and upper bounds, we recall 
the lower bound (1.7) of Harremo\"es, Johnson and Kontoyiannis \cite{H-J-K}.

\vskip5mm
{\bf Proposition 3.4.} {\sl If $\max_j p_j \leq 1/2$, then
$$
\frac{1}{4}\,\Big(\frac{\lambda_2}{\lambda}\Big)^2 \, \leq \,
D(W||Z) \, \leq \, \chi^2(W,Z) \, \leq \, 
C_\lambda\,\Big(\frac{\lambda_2}{\lambda}\Big)^2,
$$
where $C_\lambda$ depends on $\lambda \geq 0$ as an increasing
continuous function with $C_0 = 2$. In particular,
if $\lambda \leq 1/2$, then
$$
\chi^2(W,Z) \, \leq \, 15\,\Big(\frac{\lambda_2}{\lambda}\Big)^2.
$$
}

{\bf Proof.} 
Applying Lemmas 3.1-3.2, we get
$$
\lambda_2^{-2}\,e^{\lambda}\, \chi^2(W,Z)  \, \leq \,
1+\frac{1}{\lambda}(\lambda+e-1)^2  + \sum_{k=2}^\infty \frac{k!}{\lambda^k}
\Big(\frac{\lambda^k}{k!} + c_\lambda\,\frac{\lambda^{k-1}}{(k-1)!} + 
\frac{\lambda^{k-2}}{(k-2)!}\Big)^2,
$$
where $c_\lambda = \frac{e^\lambda - 1}{\lambda}$.
Expanding the squares of the brackets in this sum results in
$$
\sum_{k=2}^\infty \, \frac{k!}{\lambda^k} \, 
\Big(\frac{\lambda^{2k}}{k!^{\,2}} + 
\frac{2c_\lambda\,\lambda^{2k-1}}{k!\, (k-1)!} + 
\frac{c_\lambda^2\,\lambda^{2k-2}}{(k-1)!^{\,2}} + 
\frac{2\lambda^{2k-2}}{k!\,(k-2)!} + 
\frac{2c_\lambda\,\lambda^{2k-3}}{(k-1)!\,(k-2)!} + 
\frac{\lambda^{2k-4}}{(k-2)!^{\,2}}\Big)
$$
\bee
= \
\sum_{k=2}^\infty \frac{\lambda^k}{k!} + 
2c_\lambda\, \sum_{k=2}^\infty \frac{\lambda^{k-1}}{(k-1)!} + 
c_\lambda^2\, \sum_{k=2}^\infty k\,\frac{\lambda^{k-2}}{(k-1)!} +
2\, \sum_{k=2}^\infty \frac{\lambda^{k-2}}{(k-2)!} \\
 + \ 
2c_\lambda\, \sum_{k=2}^\infty k\,\frac{\lambda^{k-3}}{(k-2)!} + 
\sum_{k=2}^\infty k(k-1)\,\frac{\lambda^{k-4}}{(k-2)!}, 
\ene
which is the same as
\bee
3e^\lambda - 1 - \lambda + 2c_\lambda\,(e^\lambda - 1) + 
c_\lambda^2 \sum_{k=1}^\infty (k+1)\,\frac{\lambda^{k-1}}{k!}\\ 
 & & \hskip-50mm + \ 
2c_\lambda \sum_{k=0}^\infty (k+2)\,\frac{\lambda^{k-1}}{k!}  + 
\sum_{k=0}^\infty (k+1)(k+2)\,\frac{\lambda^{k-2}}{k!} \\
 & & \hskip-80mm = \
3e^\lambda - 1 - \lambda + 2c_\lambda\,(e^\lambda - 1) + 
2c_\lambda e^\lambda\,\frac{2+\lambda}{\lambda} + 
\frac{2 + 4\lambda + \lambda^2}{\lambda^2}\ e^\lambda.
\ene
Multiplying by $\lambda^2$, this gives the desired inequality
$$
\lambda^2 \lambda_2^{-2}\,\chi^2(W,Z) \, \leq \, C_\lambda \, = \,
\lambda^2 + \lambda(\lambda+e-1)^2 + B_\lambda
$$
with
\bee
B_\lambda 
 & = &
\lambda^2\, (3e^\lambda - 1 - \lambda) + 2\lambda\,(e^\lambda - 1)^2 + 
2\,(2+\lambda)\,e^\lambda\,(e^\lambda - 1) + 
(2 + 4\lambda + \lambda^2)\, e^\lambda \\
 & = &
\lambda\,(2 - \lambda - \lambda^2) - 
2\,(1 + \lambda - 2\lambda^2)\,e^\lambda + 4\,(1+\lambda)\,e^{2\lambda}.
\ene
It is easy to check that $\frac{d}{d\lambda}\,B_\lambda > 0$, so that
this function is increasing in $\lambda$, with $C_0 = B_0 = 2$.

For the range $\lambda\leq\frac{1}{2}$, the term $e-1$ appearing in the definition of $C_{\lambda}$ may be replaced with $2(\sqrt e-1)$ (according to the inequality (3.7)), which leads to the constant 
$C_{1/2} = \frac{1}{2}(\frac{1}{2}+2(\sqrt e-1))^2+ \frac{7}{8} - 2\sqrt{e} + 6\,e < 15$.
\qed

\vskip7mm
\section{{\bf \large Generating functions}}
\setcounter{equation}{0}

\vskip2mm
\noindent
The probability function $f(k) = \P\{Z=k\}$ of the 
Poisson random variable $Z \sim P_\lambda$ satisfies the equation 
$\lambda f(k-1) = kf(k)$ in integers $k \geq 1$, which immediately implies
$$
\lambda\, \E\,h(Z+1) = \E\,Zh(Z)
$$
for any function $h$ on $\Z$ (as long as the expectations exist).
This identity was emphasized by Chen \cite{Ch} who proposed to consider
an approximate equality
$$
\lambda\, \E\,h(X+1) \sim \E\,Xh(X)
$$
as a characterization of a random variable $X$ being almost
Poisson with parameter $\lambda$. This idea was inspired by
a similar approach of Charles Stein to problems of normal approximation 
on the basis of the approximate equality $\E\,h'(X) \sim \E\,Xh(X)$.

Another natural approach to the Poisson approximation is based on 
the comparison of characteristic functions. Since the random variables
$W$ and $Z$ take non-negative integer values, one may equivalently consider
the associated generating functions.

The generating function for the Poisson law $P_\lambda$ with parameter
$\lambda>0$ is given by
\be
\varphi(w) = \E\,w^Z = \sum_{k=0}^\infty \P\{Z=k\}\, w^k =
e^{\lambda (w-1)} = 
\prod_{j=1}^n\, e^{p_j (w-1)}, 
\en
which is an entire function of the complex variable $w$.
Correspondingly, the generating function for the distribution of
the random variable $W = X_1 + \dots + X_n$ in (1.1) is
\be
g(w) = \E\,w^W = \sum_{k=0}^\infty \P\{W=k\}\, w^k = 
\prod_{j=1}^n\, (q_j + p_j w),
\en
which is a polynomial of degree $n$. Hence, the difference between the 
involved probabilities may be expressed via the contour integrals by the
Cauchy formula
\be
\Delta_k=\P\{W=k\} - \P\{Z=k\} = \int_{|w| = r} 
w^{-k}\, (g(w)-\varphi(w))\,d\mu_r(w),
\en
where $\mu_r$ is the uniform probability measure on the circle $|w| = r$ 
of an arbitrary radius $r>0$.

Note that for $w = e^{it}$ with real $t$, the 
generating functions $\varphi$ and $g$ become the characteristic functions 
of $Z$ and $W$, respectively. Hence, closeness of the distributions of 
these random variables may be studied as a problem of the closeness of the 
generating functions on the unit circle.

Let us now describe first steps based on the application of the formula (4.3).
Given complex numbers $a_j,b_j$ ($1 \leq j \leq n$), we have an identity
\be
a_1 \dots a_n - b_1 \dots b_n = \sum_{j=1}^n\, (a_j - b_j) 
\prod_{l<j} b_l \prod_{l>j} a_l
\en
with the convention that $\prod_{l<j} b_l = 1$ for $j=1$ and 
$\prod_{l>j} a_l = 1$ for $j = n$. It implies
$$
\bigg|\prod_{j=1}^n a_j - \prod_{b=1}^n b_j\bigg| \, \leq \,
\sum_{j=1}^n |a_j - b_j|\prod_{l<j} |b_l| \prod_{l>j} |a_l|.
$$
According to the product representations (4.1)-(4.2) to be used in (4.3), 
one should choose here $a_j = q_j + p_j w$ and $b_j = e^{p_j (w-1)}$
with $|w| = r$. Then
\be
|a_j| \, \leq \, q_j + p_j r  \, \leq \, e^{p_j (r-1)}, \qquad
|b_j| = e^{p_j ({\rm Re}\,w-1)} \leq e^{p_j (r-1)}.
\en
Therefore
\begin{eqnarray}
|g(w)-\varphi(w)| 
 & \leq &
\sum_{j=1}^n |a_j - b_j|\prod_{l \neq j} e^{p_l (r-1)} \nonumber \\ 
 & = &
e^{\lambda (r-1)} \sum_{j=1}^n |a_j - b_j|\,e^{-p_j (r-1)}.
\end{eqnarray}
To estimate the terms in this sum, consider the function 
\be
\xi(u) = 1 + u - e^u = -u^2 \int_0^1 e^{tu}\, (1-t)\,dt, \qquad
u \in \C, 
\en
of the complex variable $u$, where the Taylor integral formula is applied in the second representation.
If ${\rm Re}\, u \leq 0$, then 
$
|u^2\, e^{tu}| = |u|^2\,\exp\{t\,{\rm Re}\, u\} \leq |u|^2, 
$
so, 
\be
|\xi(u)| \leq \frac{1}{2}\,|u|^2, \qquad {\rm Re}\, u \leq 0.
\en
In particular, for $u = p_j(w-1)$ with $w = \cos \theta + i \sin \theta$, 
we have
$$
|w-1|^2 = (\cos \theta - 1)^2 + \sin^2 \theta = 2 (1 - \cos \theta),
$$
hence $|\xi(u)| \leq p_j^2\,(1 - \cos \theta)$, and (4.6) yields
$$
|g(w)-\varphi(w)| \, \leq \,
\sum_{j=1}^n |\xi(p_j(w-1))| \, \leq \, (1 - \cos \theta) \sum_{j=1}^n p_j^2 
\, \leq \, (1 - \cos \theta)\, \lambda_2.
$$
Integrating over the unit circle in (4.3), we then arrive at the uniform bound:

\vskip5mm
{\bf Proposition 4.1.} {\sl We have
\be
\sup_{k \geq 0}\, |\,\P\{W=k\} - \P\{Z=k\}| \, \leq \, \lambda_2.
\en
}

This is a weakened variant of Le Cam's bound
$|\P\{W \in A\} - \P\{Z \in A\}| \, \leq \, \lambda_2$,
specialized to the one-point set $A = \{k\}$. In order to get a similar bound 
with arbitrary sets, or develop applications to stronger distances, we need 
sharper forms of (4.9), with the right-hand side properly depending on $k$.

\vskip7mm
\section{{\bf \large Proof of Theorem 1.3}}
\setcounter{equation}{0}

\vskip2mm
\noindent 
Applying (4.4) with $a_j = q_j + p_j w$ and $b_j = e^{p_j (w-1)}$
in (4.3), one may write this formula as
\be
\Delta_k \, = \, 
\P\{W=k\} - \P\{Z=k\} \, = \, \sum_{j=1}^n T_j(k), \qquad k=0,1,\dots,
\en
with
\be
T_j(k) \, = \,
\int_{|w|=r}w^{-k}\,(a_j-b_j)\,\prod_{l<j}b_l\prod_{l>j} a_l\ d\mu_r(w),
\en
where the integration is performed over the uniform probability
measure $\mu_r$ on the circle $|w| = r$. 
Let us write $w = r (\cos \theta + i \sin \theta)$, $|\theta| < \pi$,
and estimate $|T_j(k)|$ by inserting the absolute 
value sign inside the integral. Then, using (4.5), we get
\bee
|T_j(k)|
 & \leq & 
r^{-k} \int_{|w|=r}\,|a_j-b_j|\,\prod_{l<j} |e^{p_l\,(w - 1)}|\,
\prod_{l>j} |q_l + p_l w|\ d\mu_r(w) \\
 & \hskip-28mm = & \hskip-15mm 
r^{-k} \int_{|w|=r}\,|a_j-b_j|\,
\exp\Big\{(r\cos \theta - 1)\, \sum_{l=1}^{j-1} p_l\Big\}\,
\prod_{l = j+1}^n |q_l + p_l w|\ d\mu_r(w) \\ 
 & \hskip-28mm = & \hskip-15mm 
r^{-k} e^{(r - 1)\, \sum_{l=1}^{j-1} p_l}
\int_{|w|=r}\,|a_j-b_j|\,
\exp\Big\{-2r\sin^2\frac{\theta}{2}\, \sum_{l=1}^{j-1} p_l\Big\}\,
\prod_{l = j+1}^n |q_l + p_l w|\ d\mu_r(w).
\ene
Here, in order to estimate $|a_j - b_j|$, let us return to the function
$\xi(u)$ introduced in (4.7), which we need at the values 
$u_j = p_j(w - 1)$ with $|w| = r$.

\vskip5mm
{\bf Case 1}: $r \geq 1$.
Since ${\rm Re}\, u_j \leq p_j(r-1)$, we have, for any $t \in (0,1)$,
$$
|u_j^2\, e^{tu_j}| = |u_j|^2\,e^{t\,{\rm Re}\, u_j} \leq |u_j|^2\,e^{p_j t(r - 1)}
 \leq |u_j|^2\,e^{p_j (r - 1)}, 
$$
so, by (4.7), 
$$
|a_j - b_j| =|\xi(u_j)| \, \leq \, \frac{1}{2}\,p_j^2\,|w-1|^2\,\,e^{p_j (r - 1)}.
$$

\vskip2mm
{\bf Case 2}: $0 < r < 1$. Then
${\rm Re}\, u_j \leq 0$, so, by (4.8),
$$
|a_j - b_j| =|\xi(u_j)| \, \leq \, \frac{1}{2}\,p_j^2\,|w-1|^2.
$$

\vskip2mm
Since $|w-1|^2 = (r-1)^2 + 4r\,\sin^2(\theta/2)$, we therefore obtain from (5.2)
that
\be
|T_j(k)| \, \le \, \frac{1}{2} \,p_j^2\,R_j(r)\,r^{-k}\, 
\Big((r-1)^2\,I_{j0}(r) + 4r I_{j2}(r)\Big),
\en
where 
$$
R_{j}(r) =
\left\{
\begin{array}{lc}
\exp\big\{(r-1) \sum_{l=1}^j p_l\big\}\, \prod_{l=j+1}^n (q_l + p_l r)
 & \ \ {\rm for} \ \  r \geq 1, \\
\exp\big\{(r-1) \sum_{l=1}^{j-1} p_l\big\}\, \prod_{l=j+1}^n (q_l + p_l r)
 & \ \ {\rm for} \ \ r < 1,
\end{array}
\right.
$$
and
$$
I_{jm}(r) \, = \,
\frac{1}{2\pi} \int_{-\pi}^{\pi} \Big|\sin\frac{\theta}{2}\Big|^m\,
\exp\Big\{-2 r\sin^2\frac{\theta}{2}\, \sum_{l=1}^{j-1} p_l\Big\}
\prod_{l=j+1}^n \frac{|q_l + p_l r\, e^{i\theta}|}{q_l+p_lr}\ d\theta.
$$

In order to estimate the last integrals, which we need
with $m=0$ and $m=2$, let us first note that 
\bee
|q_l + p_l r e^{i\theta}|^2
 & = &
q_l^2 + p_l^2 r^2 + 2p_lq_l\, r\cos\theta
 \, = \, 
(q_l + p_l r)^2 - 4 q_l p_l\, r\sin^2\frac{\theta}{2}\\
 & = &
(q_l+p_l r)^2 
\Big(1 - \frac{4q_l p_l r}{(q_l + p_l r)^2}\, \sin^2\frac{\theta}{2}\Big).
\ene
Hence, using $1 - x \leq e^{-x}$ ($x \in \R$), we have 
\begin{eqnarray}
\prod_{l=j+1}^n \frac{|q_l + p_l r\, e^{i\theta}|}{q_l+p_lr}
 & = &
\prod_{l=j+1}^n 
\Big(1-\frac{4q_l p_l r}{(q_l+p_lr)^2} \sin^2\frac{\theta}2\Big)^{1/2} \nonumber \\
 & \le &
\exp\Big\{-2\,\sin^2\frac{\theta}2\,\sum_{l=j+1}^n
\frac{q_l p_l\, r}{(q_l + p_l r)^2}\Big\},
\end{eqnarray}
so that
\begin{eqnarray}
I_{jm}(r) 
 & \leq &
\frac{1}{2\pi} \int_{-\pi}^{\pi} \Big|\sin\frac{\theta}{2}\Big|^m\,
\exp\Big\{-2\gamma_j(r)\,\sin^2\frac{\theta}{2}\Big\}\ d\theta \nonumber \\
 & \leq &
\frac{1}{2\pi}\, 2^{-m} \int_{-\pi}^{\pi} |\theta|^m
\exp\Big\{-\frac{2}{\pi^2}\, \gamma_j(r)\,\theta^2\Big\}\,d\theta.
\end{eqnarray}
Here we applied the inequalities
$\frac{2}{\pi}\,t \leq \sin t \le t$ ($0\le t\le\frac{\pi}{2}$) and used
the notation
$$
\gamma_j(r) = r\,\Big(\sum_{l=1}^{j-1} p_l + 
\sum_{l=j+1}^n \frac{q_l p_l}{(q_l + p_l r)^2}\Big).
$$

Thus, we need to bound $\gamma_j$ from below. If $r \geq 1$, then 
$q_l + p_l r \leq r$, so
$$
\sum_{l=j+1}^n \frac{q_l p_l\, r}{(q_l + p_l r)^2}\ge 
\frac{1}{r} \sum_{l=j+1}^n q_l p_l.
$$
This gives
\bee
\gamma_j(r) 
 & \geq &
r \sum_{l=1}^{j-1} p_l + \frac{1}{r} \sum_{l=j+1}^n q_l p_l \\ 
 & = &
r \sum_{l=1}^{j-1} p_l + \frac{1}{r} \sum_{l=1}^n (p_l - p_l^2) -
\frac{1}{r} \sum_{l=1}^j (p_l - p_l^2) \\ 
 & = &
\Big(r - \frac{1}{r}\Big) \sum_{l=1}^{j-1} p_l + 
\frac{1}{r} \sum_{l=1}^{j-1} p_l^2 + \frac{1}{r} \sum_{l=1}^n (p_l - p_l^2)
- \frac{1}{r}\,q_j p_j
 \, \geq \,
\frac{1}{r}\,(\lambda - \lambda_2 - q_j p_j).
\ene
In case $r\le 1$, we use $q_l+p_lr\le 1$, implying that
$$
\sum_{l=j+1}^n \frac{q_l p_l}{(q_l + p_l r)^2} \ge \sum_{l=j+1}^n q_l p_l.
$$
Therefore in this range we have a similar lower bound, namely
\bee
\gamma_j(r) 
 & \geq &
r \sum_{l=1}^{j-1} p_l + r \sum_{l=j+1}^n q_l p_l \\ 
 & = &
r \sum_{l=1}^{j-1} p_l + r \sum_{l=1}^n (p_l - p_l^2) -
r \sum_{l=1}^j (p_l - p_l^2) \\ 
 & = &
-r p_j + r \sum_{l=1}^{j} p_l^2 + r \sum_{l=1}^n (p_l - p_l^2)
 \, \geq \,
r\,(\lambda - \lambda_2 - q_j p_j).
\ene
Since $q_j p_j \leq \frac{1}{4}$, both lower bounds yield
$$
\gamma_j(r) \, \geq \, \psi(r) - \frac{1}{4}, \qquad
\psi(r) = \min\{r,1/r\}\,(\lambda - \lambda_2).
$$

As a result, (5.5) is simplified to
\bee
I_{jm}(r) 
 & \leq &
\frac{1}{2\pi}\, 2^{-m} \sqrt{e} \int_{-\pi}^{\pi} |\theta|^m
\exp\Big\{-\frac{2}{\pi^2}\, \psi(r)\,\theta^2\Big\}\,d\theta \\
 & = &
\sqrt{e}\,\frac{\pi^m}{4^{m+1}}\,\psi(r)^{-\frac{m+1}{2}}
\int_{-2\sqrt{\psi(r)}}^{2\sqrt{\psi(r)}}\,|x|^m\, e^{-\frac{1}{2}\,x^2}\,dx.
\ene
The last integral may be extended to the whole real line, which makes sense 
for large values of $\psi(r)$, or one may bound the exponential term in 
the integrand by 1, which makes sense for small values of $\psi(r)$. 
These two ways of estimation lead to
\bee
I_{jm}(r) 
 & \leq &
\sqrt{e}\,\frac{\pi^m}{4^{m+1}}\,\psi(r)^{-\frac{m+1}{2}}
\min\Big\{\sqrt{2\pi}\,
\E\,|\xi|^m, \frac{2^{m+2}}{m+1}\,\psi(r)^{\frac{m+1}{2}}\Big\} \\
 & \leq &
\sqrt{e}\,\frac{\pi^m}{4^{m+1}}\,
\max\Big\{\sqrt{2\pi}\,\E\,|\xi|^m, \frac{2^{m+2}}{m+1}\Big\}\,
\min\Big\{1,\psi(r)^{-\frac{m+1}{2}}\Big\},
\ene
where $\xi$ is a standard normal random variable. In particular,
we get the upper bounds
$$
I_{j0}(r) \, \le \,
\sqrt{e}\, \min\big\{1,\psi(r)^{-1/2}\big\}, \qquad
I_{j2}(r) \, \le \, 
\frac{\sqrt{e}\,\pi^2}{12}\, \min\big\{1,\psi(r)^{-3/2}\big\}.
$$

In view of $q_l + p_l r \leq e^{(r-1) p_l}$, from the definition of $R_j(r)$
we also have the bound
$$
R_j(r) \, \le \,
\exp\Big\{(r-1)\sum_{l=1}^n p_l\Big\} \, = \, e^{\lambda(r-1)}
$$
in case $r \geq 1$, while for $r \leq 1$
$$
R_j(r) \, \le \,
\exp\Big\{(r-1)\sum_{l \neq j} p_l\Big\} \, = \, e^{\lambda(r-1)}\,
e^{-p_j (r-1)} \, \leq \, e^{\lambda(r-1) + 1}.
$$
Applying these bounds in (5.3), we therefore obtain that $|T_j(k)|$
may be bounded from above by
$$
\frac{\delta_r}{2}\, p_j^2\,e^{\lambda(r-1) + \frac{1}{2}}\,r^{-k}
\Big((r-1)^2 \min\big\{1,\psi(r)^{-1/2}\big\} + 
\frac{\pi^2}{3}\, r \min\big\{1,\psi(r)^{-3/2}\big\}\Big),
$$
where $\delta_r = 1$ in case $r \geq 1$ and $\delta_r = e$ for $r < 1$.
Summing over $j \leq n$ and recalling (5.1), one can estimate 
$|\Delta_k|$ from above by
\be
\lambda_2\, \delta_r\,e^{\lambda(r-1)}\,r^{-k}\,
\Big(\frac{\sqrt{e}}{2}\, (r-1)^2\, \min\big\{1,\psi(r)^{-1/2}\big\} + 
\frac{\sqrt{e}\,\pi^2}{6}\ r\,\min\big\{1,\psi(r)^{-3/2}\big\}\Big).
\en

Now, letting $r \rightarrow 0$ in the case  
$k=0$, (5.6) leads to 
$$
|\Delta_0|\le \frac{e\sqrt e}{2}\,\lambda_2\,e^{-\lambda}<3\lambda_2\,e^{-\lambda},
$$
and we obtain the first inequality in (1.9).
Letting $r\downarrow 1$ in the case $k\geq 1$, (5.6)
gives
$$
|\Delta_k| \leq \frac{\sqrt{e}\,\pi^2}{6} \,\lambda_2 < 
3\lambda_2,
$$
which is the second inequality in (1.9). 

But, if $k\geq 1$, one may also use (5.6) with $r = \frac{k}{\lambda}$ and apply the bound
$k! \leq e\,k^{k + \frac{1}{2}}\,e^{-k}$, cf. (2.1), giving
$$
e^{\lambda(r-1)}\,r^{-k} = 
\Big(\frac{e \lambda}{k}\Big)^k\,e^{-\lambda} \leq e\sqrt{k}\,f(k), \qquad
f(k) = \frac{\lambda^k}{k!}\,e^{-\lambda}.
\qquad 
$$
To simplify the numerical constants, note that $\frac{1}{2}\,e^{5/2} < 6.1$ 
and $\frac{1}{6}\,e^{5/2}\,\pi^2 < 20.1$. Recalling that
$\psi(r) = \rho$ for $r = k/\lambda$, we finally get
the second inequality (1.10),
\be
|\Delta_k| \, \leq \,
\lambda_2 \sqrt{k}\,f(k)\,
\Big(7\,\Big(\frac{k - \lambda}{\lambda}\Big)^2\, \min\big\{1,\rho^{-1/2}\big\} + 
21\,\frac{k}{\lambda}\,\min\big\{1,\rho^{-3/2}\big\}\Big).
\en
\qed

\vskip7mm
\section{{\bf \large Consequences of Theorem 1.3}}
\setcounter{equation}{0}

\vskip2mm
\noindent 
Under the natural requirement that $\lambda_2$ is bounded away from 
$\lambda$, the bound (5.7) on $\Delta_k = \P\{W=k\}-\P\{Z=k\}$
may be simplified. As before, we use the notations
$$
f(k) = \P\{Z=k\} = \frac{\lambda^k}{k!}\,e^{-\lambda}, \qquad
\lambda = p_1 + \dots + p_n, \quad \lambda_2 = p_1^2 + \dots + p_n^2.
$$
Note that $\lambda_2 \leq \lambda$ and recall that
$\rho = (\lambda-\lambda_2)\,\min\{\frac{k}{\lambda},\frac{\lambda}{k}\}$.

\vskip5mm
{\bf Corollary 6.1.} {\sl If $\lambda_2 \leq \kappa \lambda$,
$\kappa \in (0,1)$, then for any integer $k \geq 0$,
\be
\hskip-0mm
|\Delta_k| \, \leq \,
\frac{7 f(k)}{(1 - \kappa)^{3/2}}\,
\Big(\frac{(k-\lambda)^2}{\lambda} + 3\Big)\, \frac{\lambda_2}{\lambda}\,
\max\Big\{\Big(\frac{k}{\lambda}\Big)^3,1\Big\}.
\en
In particular, if $k \leq 2\lambda$, then
\be
|\Delta_k| \, \leq \,
\frac{56\, f(k)}{(1 - \kappa)^{3/2}}\,\Big(\frac{(k-\lambda)^2}{\lambda} + 3\Big)\, 
\frac{\lambda_2}{\lambda}.
\en
If $k \geq \lambda \geq 1/2$, we also have
\be
|\Delta_k| \, \leq \,
\frac{49\, f(k)}{(1 - \kappa)^{3/2}}\,\Big(\frac{k}{\lambda}\Big)^3\,\lambda_2.
\en
}

\vskip2mm
{\bf Proof.}
The assumption $\lambda_2 \leq \kappa \lambda$ ensures that
$\rho \geq (1 - \kappa)\lambda\,\min\{\frac{k}{\lambda},\frac{\lambda}{k}\}$.

If $1 \leq k \leq K \lambda$ ($K \geq 1$), then
$\frac{k}{\lambda} \leq K^2 \frac{\lambda}{k}$ and
$\rho \geq \frac{1 - \kappa}{K^2}\,k$, so, the right-hand side of
(5.7) is bounded from above by
$$
\lambda_2 \sqrt{k}\,f(k)\,
\Big(7\,\Big(\frac{k - \lambda}{\lambda}\Big)^2\, 
\frac{K}{\sqrt{(1 - \kappa)\,k}} + 
21\,\frac{k}{\lambda}\,\frac{K^3}{(1 - \kappa)^{3/2}\,k^{3/2}}\Big).
$$
Choosing $K = \max\{\frac{k}{\lambda},1\}$, this expression does not exceed 
the right-hand side of (6.1). Thus, the inequality (1.10) yields (6.1), 
which in turn immediately implies (6.2).

In case $k = 0$, we apply the inequality (1.9). Since
$\frac{(k-\lambda)^2}{\lambda} + 3 \geq \lambda$ for $k = 0$,
the right-hand side of (1.10) is dominated by the right-hand side of (6.1).
Thus, we obtain (6.1) without any constraints on $k$, and (6.2) for all 
$k \leq 2 \lambda$.

In case $k \geq \lambda$, necessarily $\rho \geq (1 - \kappa)\,\lambda^2/k$. 
Hence, the right-hand side of (5.7) may be bounded from above by
$$
\lambda_2 \sqrt{k}\,f(k)\,
\Big(7\,\Big(\frac{k - \lambda}{\lambda}\Big)^2\, 
\frac{\sqrt{k}}{\lambda \sqrt{1 - \kappa}} + 
21\,\frac{k}{\lambda}\cdot \frac{k^{3/2}}{\lambda^3\,(1 - \kappa)^{3/2}}\Big).
$$
Using $(\frac{k-\lambda}{\lambda})^2 \leq \frac{k^2}{\lambda^2}$ to bound
the first term in the brackets and $\frac{k}{\lambda} \leq 2k$ to bound 
the second term (using $\lambda \geq 1/2$), we obtain the bound (6.3).
\qed

\vskip5mm
We are now prepared to extend Proposition 3.4 to larger values
of $\lambda$ under the assumption that $\lambda_2/\lambda$ is
bounded away from 1. The next assertion, being combined with 
Proposition 3.4, yields Theorem 1.2 with $c = 15$ in case 
$\lambda \leq 1/2$ and $c = 56 \cdot 10^6$ in case $\lambda > 1/2$.

\vskip5mm
{\bf Proposition 6.2.} {\sl If $\lambda \geq 1/2$ and
$\lambda_2 \le \kappa \lambda$ with $\kappa \in (0,1)$, then
\be
\frac{1}{4}\,\Big(\frac{\lambda_2}{\lambda}\Big)^2 \, \leq \, 
 D(W||Z) \, \leq \, \chi^2(W,Z) \, \leq \,
c_\kappa\,\Big(\frac{\lambda_2}{\lambda}\Big)^2.
\en
where $c_\kappa = c\,(1-\kappa)^{-3}$ with, for example, $c = 7 \cdot 10^6$.
}

\vskip5mm
{\bf Proof.} 
The leftmost lower bound in (6.4) is added according to (1.7)
(using the Pinsker inequality, it also follows with some
constant from Barbour-Hall's lower bound in Theorem 1.1). Hence, it remains to show
the rightmost upper bound in (6.4). Write
$$
\chi^2(W,Z) \, = \, \sum_{k=0}^{\infty} \frac{\Delta_k^2}{f(k)} \, = \, 
S_1 + S_2 \, = \,
\bigg(\,\sum_{k=0}^{[2\lambda]} + \sum_{k=[2\lambda]+1}^{\infty}\bigg)
\frac{\Delta_k^2}{f(k)}.
$$

In the range $0\le k\le [2\lambda]$, we apply the inequality (6.2) 
which gives
$$
\Delta_k^2 \, \leq \, \frac{56^2}{(1 - \kappa)^3}\,
\Big(\frac{(k - \lambda)^4}{\lambda^2} + 
6\,\frac{(k - \lambda)^2}{\lambda^2} + 9\Big)\,
\Big(\frac{\lambda_2}{\lambda}\Big)^2\,f(k)^2.
$$
Hence
$$
S_1
 \, \le \, \frac{56^2}{(1 - \kappa)^3}\,
\Big(\frac{\E\,(Z - \lambda)^4}{\lambda^2} + 6\,
\frac{\E\,(Z - \lambda)^2}{\lambda} + 9\Big)\,
\Big(\frac{\lambda_2}{\lambda}\Big)^2.
$$
In the sequel, we use a simple moment inequality
$\E\,Z^m \leq \lambda (\lambda + 1) \dots (\lambda + m - 1)$.
We also have $\E\,(Z - \lambda)^2 = \lambda$ and
$\E\,(Z - \lambda)^4 = \lambda (\lambda+3)$, so that
\begin{eqnarray}
S_1
 & \le & 
\frac{56^2}{(1 - \kappa)^3}\,
\Big(\frac{\lambda + 3}{\lambda} + 15\Big)\,
\Big(\frac{\lambda_2}{\lambda}\Big)^2 \nonumber \\
 & \leq &
\frac{18\,816}{(1 - \kappa)^3}\,\big(\lambda^{-1} + 3\big)\,
\Big(\frac{\lambda_2}{\lambda}\Big)^2 \, \leq \, 
\frac{C_1}{(1 - \kappa)^3}\,\Big(\frac{\lambda_2}{\lambda}\Big)^2
\end{eqnarray}
with $C_1 = 94\,080$
(where we used the assumption $\lambda \geq 1/2$ on the last step).

In order to estimate $S_2$, we use the following elementary bound
\be
\sum_{k = k_0}^\infty k^d f(k) \, \leq \, k_0^d\,f(k_0)\,
\Big(1 - \frac{\lambda}{k_0}\,\Big(\frac{k_0 + 1}{k_0}\Big)^{d-1}\Big)^{-1},
\en
which holds for any $d=1,2,\dots$ as long as 
$k_0^d/(k_0 + 1)^{d-1} > \lambda$. For the proof, write
$$
\sum_{k = k_0}^\infty k^d f(k) \, = \, k_0^d f(k_0)\,
\big(1 + \theta_1 + \theta_1 \theta_2 + \dots + 
\theta_1 \dots \theta_m + \dots\big),
$$
where 
$$
\theta_m = 
\Big(\frac{k_0 + m}{k_0 + m - 1}\Big)^d\,\frac{\lambda}{k_0 + m}, 
 \qquad m = 1,2,\dots
$$ 
Since the function $(x+1)^{d-1}\,x^{-d}$ is 
decreasing in $x>0$, we have
$1 > \theta_1 > \theta_2 > \dots$ This gives
$$
\sum_{k = k_0}^\infty k^d f(k) \, \leq \, k_0^d f(k_0)\,
\Big(1 + \sum_{m=1}^\infty \theta_1^m\Big),
$$
that is, (6.6).
In particular, for $k_0 = [2\lambda] + 1$ and $\lambda \geq 8$ (with $d=6$),
$$
\Big(1 - \frac{\lambda}{k_0}\,\Big(\frac{k_0 + 1}{k_0}\Big)^5\Big)^{-1} < 
\Big(1 - \frac{1}{2}\,\Big(\frac{2\lambda + 1}{2\lambda}\Big)^5\Big)^{-1} <
3.1.
$$
So, by (6.6), and using $[2\lambda] + 1 \leq \frac{17}{8}\lambda$
for the chosen range of $\lambda$, we have
$$
\sum_{k = [2\lambda] + 1}^\infty k^6 f(k) \, \leq \, 
3.1\,([2\lambda] + 1)^6\,f([2\lambda] + 1) \, \leq \, 
3.1 \cdot (17\lambda/8)^6\,f([2\lambda] + 1).
$$
Hence, by (6.3), 
\be
S_2 \, = \,
\sum_{k=[2\lambda]+1}^{\infty} \frac{|\Delta_k|^2}{f(k)}\, \leq \,
\frac{49^2}{(1 - \kappa)^3}\sum_{[2\lambda]+1}^{\infty}
\Big(\frac{k}{\lambda}\Big)^6\,  \lambda_2^2\, f(k)
 \, \leq \,
\frac{C_2\,\lambda_2^2}{(1 - \kappa)^3}\, f([2\lambda] + 1)
\en
with $C_2 = 49^2 \cdot 3.1 \cdot (17/8)^6 < 685\,343$.
Asymptotically with respect to large $\lambda$, this bound is much better
than (6.4). Applying $f(k) \le \frac{1}{\sqrt{2\pi k}}\, 
e^{k - \lambda}\, (\frac{\lambda}{k})^k$ as in (2.5)
with $k = [2\lambda] + 1$ and using $2\lambda \leq k \leq 2\lambda+1$, 
we have
$$
f([2\lambda] + 1) \, \leq \,
\frac{e}{2\sqrt{\lambda \pi}}\, (e/4)^\lambda \, \leq \,
\frac{e}{2\sqrt{\pi}}\ 8^{3/2}\,\Big(\frac{e}{4}\Big)^8\,
\frac{1}{\lambda^2} \, < \, \frac{1}{\lambda^2}.
$$
This gives
$$
S_2 \leq \frac{C_2}{(1 - \kappa)^3}\,\Big(\frac{\lambda_2}{\lambda}\Big)^2.
$$
As a result, we arrive at the desired upper bound in (6.4).

Finally, let us estimate $S_2$ for the range 
$\frac{1}{2} \leq \lambda \leq 8$. Returning to (6.7), we have 
$$
S_2 \, \leq \, \frac{49^2}{(1 - \kappa)^3}\sum_{k=1}^{\infty}
\Big(\frac{k}{\lambda}\Big)^6\, \lambda_2^2\, f(k) \, \leq \,
\frac{49^2}{(1 - \kappa)^3}\,\lambda^{-6}\,\lambda_2^2\,\E Z^6  \, \leq \,
\frac{C_2'}{(1 - \kappa)^3}\,\Big(\frac{\lambda_2}{\lambda}\Big)^2,
$$
where 
$C_2' = 49^2\,\sup_{\frac{1}{2} \leq \lambda \leq 8} \psi(\lambda)$,
$\psi(\lambda) = \lambda^{-4}\,\E Z^6$. Here
$$
\psi(\lambda) = \frac{(\lambda+1)\dots(\lambda+5)}{\lambda^3} =
\psi_1(\lambda) \psi_2(\lambda) \psi_3(\lambda)
$$
with $\psi_1(\lambda) = 5 + \lambda + \frac{4}{\lambda}$, 
$\psi_2(\lambda) = 7 + \lambda + \frac{10}{\lambda}$, 
$\psi_3(\lambda) = 1 + \frac{3}{\lambda}$. All these three functions are
convex, while $\psi_3$ is decreasing. In addition, $\psi_i(1/2) \geq \psi_i(8)$
for $i = 1,2$. Hence $\psi(\lambda) \leq \psi(1/2) = \frac{1}{4}\cdot 11!!$
It follows that
$C_2' = 49^2\cdot \frac{1}{4}\cdot 11!! < 6\,239\,560$, and thus
$c = C_1 + C_2''$ is the resulting constant in (6.4).
\qed

\vskip5mm
{\bf Remark 6.3.} Up to a numerical constant, the upper bound in (6.4)
immediately implies an upper bound of Theorem 1.1 in case $\lambda \geq 1/2$, 
in view of the relation $d(W,Z)^2 \leq \frac{1}{2}\,D(W,Z)$.
Indeed, (6.4) gives $d(W,Z) \leq c_\kappa \lambda_2/\lambda$,
provided that $\lambda_2 \le \kappa \lambda$. But, in the other case 
$\lambda_2 \ge \kappa \lambda$, there is nothing to prove, since
$d(W,Z) \leq 2$. Note also that, for $\lambda \leq 1/2$, 
the correct upper bound on the total variation distance is of the form 
$d(W,Z) \leq C\lambda_2$. It may be obtained as  
a consequence of Lemmas 3.1-3.2.

\vskip7mm
\section{{\bf \large Uniform Bounds. Comparison with Normal Approximation}}
\setcounter{equation}{0}

\vskip2mm
\noindent
A different choice of the parameter $r$ in the proof of Theorem 1.3 may provide 
various uniform bounds in the Poisson approximation, like in the next 
assertion. Using the $L^\infty(\mu)$-norm with respect to the counting measure 
$\mu$ on $\Z$, let us focus on the deviations of the densities of $W$ and $Z$
and the deviations of their distribution functions. These distances are
thus given by
\bee
M(W,Z) 
 & = &
\sup_{k \geq 0}\ |\P\{W = k\} - \P\{Z = k\}|, \\ 
K(W,Z)
 & = &
\sup_{k \geq 0}\ |\P\{W \leq k\} - \P\{Z \leq k\}|.
\ene
Putting $r = 1$ in (5.6), we arrive at the next assertion 
which sharpens Proposition~4.1.

\vskip5mm
{\bf Theorem 7.1.} {\sl We have
\be
M(W,Z) \, \leq \, \frac{\sqrt{e}\,\pi^2}{6}\, \lambda_2 
\min\big\{1,(\lambda - \lambda_2)^{-3/2}\big\}.
\en
}

\vskip2mm
This uniform bound is not new; with a non-explicit numerical factor, it 
corresponds to Theorem~3.1 in Cekanavicius \cite{Ce}, p.\,53.
For $\lambda \leq 1$, this relation is simplified to
$$
M(W,Z) \, \leq \, \frac{\sqrt{e}\,\pi^2}{6}\, \lambda_2,
$$
which cannot be improved (modulo a numerical factor) in view of 
the lower bounds on $|\Delta_k|$ with $k = 0,1,2$ mentioned in Section 3.
We also have a similar bound for the Kolmogorov distance,
$K(W,Z) \leq C\lambda_2$, which follows from the upper bound for 
the stronger total variation distance as in Theorem 1.1.

When, however, $\lambda$ is large (and say all $p_j \leq 1/2$), one would expect to achieve more accurate bounds when replacing the Poisson
approximation for $P_W$ by the normal law $N(\lambda,\lambda)$ with mean 
$\lambda$ and variance $\lambda$. Indeed, suppose, for example, that 
$p_j = 1/2$, so that $W$ has a binomial distribution with parameters
$(n,1/2)$, while the approximating Poisson distribution has parameter 
$\lambda = n/2$ with $\lambda_2 = n/4$. Here (1.2) only yields
$d(W,Z) \sim 1$, which means that there is no Poisson approximation
with respect to the total variation! Nevertheless, the approximation is still 
meaningful in a weaker sense in terms of the Kolmogorov distance $K$, 
as well as in terms of $M$. In this case, both $P_W$ and $P_\lambda$ 
are almost equal to $N(\lambda,\lambda)$, and the Berry-Esseen theorem 
provides a correct bound $K(W,Z) \leq \frac{c}{\sqrt{n}}$ via the triangle 
inequality for~$K$. Since $M \leq 2K$ (which holds true for all probability 
distributions on $\Z$), we also have $M(W,Z) \leq \frac{c}{\sqrt{n}}$.
Note that this inequality also follows from Theorem 7.1.
Indeed, when $\lambda_2 \leq \frac{1}{2}\,\lambda$, (7.1) is simplified to
\be
M(W,Z) \, \leq \, \frac{\sqrt{2e}\,\pi^2}{3}\, \frac{\lambda_2}{\lambda^{3/2}},
\en
which yields a correct order for growing $n$. Thus, the two approaches
are equivalent for this particular (i.i.d.) example.

To realize whether or not the normal approximation is better or worse than 
the Poisson approximation in the general non-i.i.d. situation (that is, with different
$p_j$'s), let us evaluate the corresponding Lyapunov ratio in the central 
limit theorem and apply the Berry-Esseen bound
$K(W,N_\lambda) \leq cL_3$, where the random variable
$N_\lambda$ is distributed according to $N(\lambda,\lambda)$.
Since $\Var(W) = \sum_{j=1}^n p_j q_j = \lambda - \lambda_2$, the Lyapunov
ratio for the sequence $X_1,\dots, X_n$ is given by
\bee
L_3 
 & = &
\frac{1}{\Var(W)^{3/2}}\,\sum_{j=1}^n \E\,|X_j - \E X_j|^3 \\
 & = &
\frac{1}{(\lambda - \lambda_2)^{3/2}}\,\sum_{j=1}^n\, (p_j^2 + q_j^2)\,p_j q_j
 \ \leq \ \frac{1}{\sqrt{\lambda - \lambda_2}}
\ene
(note that $\frac{1}{2} \leq p_j^2 + q_j^2 \leq 1$). Hence
$K(W,N_\lambda) \leq \frac{c}{\sqrt{\lambda - \lambda_2}}$,
up to some absolute constant $c>0$. A similar bound holds for $Z$ as well
when representing $W$ as the sum of $n$ independent Poisson random variables
$Z_j$ with parameters $p_j$. Namely, for the sequence $Z_1,\dots,Z_n$, we have
$$
L_3 \, = \, \frac{1}{\Var(Z)^{3/2}}\,\sum_{j=1}^n \E\,|Z_j - \E Z_j|^3 \, \leq \,
\frac{c}{\lambda^{3/2}}\,\sum_{j=1}^n p_j \, = \, \frac{c}{\sqrt{\lambda}}.
$$
Therefore, $K(Z,N_\lambda) \leq \frac{c}{\sqrt{\lambda}}$ and hence, by the
triangle inequality, $K(W,Z) \leq \frac{c}{\sqrt{\lambda - \lambda_2}}$.
In particular, in a typical situation where $\lambda_2 \leq \frac{1}{2}\,\lambda$, 
the normal approximation yields
\be
M(W,Z) \, \leq \, \frac{c}{\sqrt{\lambda}}
\en 
with some absolute constant $c$. But, this bound is surprisingly 
worse than (7.2) as long as $\lambda_2 = o(\lambda)$.

Consider as an example $p_j = 1/(2\sqrt{j})$ for $j = 1,\dots,n$. 
Then $\lambda \sim \sqrt{n}$, $\lambda_2 \sim \log n$, and we get
$M(W,Z) \leq c n^{-3/4} \log n$ in (7.2), while (7.3) only yields
$M(W,Z) \leq c n^{-1/4}$. This example is also illustrative when
comparing Theorem 1.2 with (1.5). The first one provides
a correct asymptotic $D(W,Z) \sim \frac{\log^2 n}{n}$ (within absolute
factors), while (1.5) only gives $D(W,Z) \leq c$.

\vskip5mm
{\bf Acknowledgement.} The authors would like to thank Igal Sason
and two referees for valuable comments and
drawaing our attention to additional references related to 
the Poisson approximation in informational distances.


\vskip10mm
\begin{thebibliography}{BH3}
\itemsep=-0pt
\small
\vskip2mm 
\bibitem{B-H} 
Barbour, A. D.; Hall, P. On the rate of Poisson convergence. 
        Math. Proc. Cambridge Philos. Soc. 95 (1984), no. 3, 473--480. 

\vskip2mm 
\bibitem{B-H-J} 
Barbour, A. D.; Holst, L.; Janson, S. Poisson approximation. 
        Oxford Studies in Probability, 2. Oxford Science Publications. 
				The Clarendon Press, Oxford University Press, New York, 1992. x+277 pp. 

\vskip2mm 
\bibitem{Ce}
\v Cekanavicius, V. Approximation methods in Probability Theory.
        Universitext. Springer (2016), 274 pp.

\vskip2mm 
\bibitem{Ch}
Chen, L. H. Y. Poisson approximation for dependent trials. 
        Ann. Probability 3 (1975), no. 3, 534--545. 

\vskip2mm       
\bibitem{E-H} 
van Erven, T., Harremo\"es, P. R\'enyi divergence and Kullback-Leibler divergence.
        IEEE Trans. Inform. Theory 60 (2014), no. 7, 3797--3820.

\vskip2mm       
\bibitem{H2001} 
Harremo\"es, P. Binomial and Poisson distributions as maximum 
        entropy distributions. IEEE Trans. Inform. Theory 47 (2001), 
        no. 5, 2039--2041.

\vskip2mm       
\bibitem{H2003} 
Harremo\"es, P. Convergence to Poisson distribution in information divergence. 
        Preprint 2, Math. Department, University of Copenhagen, Feb. 2003.

\vskip2mm       
\bibitem{H-J-K} 
Harremo\"es, P.; Johnson, O.; Kontoyiannis. Thinning and information
        projections. arXive:1601.04255, Jan. 2016.

\vskip2mm       
\bibitem{H-R} 
Harremo\"es, P.; Ruzankin, P. S. Rate of convergence to Poisson law 
        in terms of information divergence. IEEE Trans. Inform. Theory 
			 50 (2004), no. 9, 2145--2149.

\vskip2mm 
\bibitem{H-LC}
Hodges, J. L., Jr.; Le Cam, L. The Poisson approximation to the Poisson 
        binomial distribution. Ann. Math. Statist. 31 (1960), 737--740. 

\vskip2mm 
\bibitem{J}        
Johnson, O. Information theory and central limit theorem. 
        Imperial College Press, London, 2004, 209 p. 

\vskip2mm 
\bibitem{K} 
Kerstan, J. Verallgemeinerung eines Satzes von Prochorow und Le Cam. 
        (German) Z. Wahrscheinlichkeitstheorie und Verw. Gebiete 
				2 (1964) 173--179. 

\vskip2mm 
\bibitem{K-H-J} 
Kontoyiannis, I.; Harremo\"es, P.; Johnson, O. Entropy and the law of 
        small numbers. IEEE Trans. Inform. Theory 51 (2005), no. 2, 466--472.

\vskip2mm 
\bibitem{LC1} 
Le Cam, L. An approximation theorem for the Poisson binomial distribution. 
        Pacific J. Math. 10 (1960), 1181--1197.

\vskip2mm
\bibitem{LC2}
LeCam, L. M.  Asymptotic Methods in Statistical Decision Theory. Springer Series 
        in Statistics. Springer-Verlag, New York, 1986. xxvi+742 pp.

\vskip2mm
\bibitem{L-V} 
Liese, F.; Vajda, I. Convex statistical distances. With German, French and Russian 
        summaries. Teubner-Texte zur Mathematik [Teubner Texts in Mathematics], 
        95. BSB B. G. Teubner Verlagsgesellschaft, Leipzig, 1987, 224 pp.

\vskip2mm 
\bibitem{P}
Prokhorov, Yu. V. Asymptotic behavior of the binomial distribution. 
        (Russian) Uspehi Matem. Nauk (N.S.) 8 (1953), no. 3 (55), 
				135--142. 

\vskip2mm 
\bibitem{R}
Romanowska, M. A note on the upper bound for the distrance in total 
        variation between the binomial and the Poisson distribution. 
				Statistica Neerlandica 31 (1977), no. 3, 127--130. 

\vskip2mm 
\bibitem{Ro}
Roos, B. On the rate of multivariate Poisson convergence. J. Multivariate Anal.  
         69  (1999),  no. 1, 120--134. 

\vskip2mm 
\bibitem{S1}
Sason, I. Improved lower bounds on the total variation distance for 
         the Poisson approximation. Statist. Probab. Lett.  83  (2013),  no. 10, 
         2422--2431.

\vskip2mm 
\bibitem{Se}
Serfling, R. J. A general Poisson approximation theorem. Ann. Probab.
        3 (1975), no. 4, 726--731. 

\vskip2mm
\bibitem{Va} 
Vajda, I. Theory of Statistical Inference and Information. Kluwer Academic 
        Publishers, Dordrecht-Borston-London, 1989, 432 pp.

\vskip2mm 
\bibitem{Ve}
Vervaat, W. Upper bounds for the distance in total variation between 
        the binomial or negative binomial and the Poisson distribution. 
				Statistica Neerlandica 23 (1969), 79--86.




    



\end{thebibliography}
\end{document}